\documentclass{scrartcl}

\usepackage{color} 
\usepackage{graphicx} 
\usepackage{latexsym,amssymb,amsmath,amsfonts} 
\usepackage{mathrsfs} 
\usepackage{bbm} 
\usepackage{todonotes}
\usepackage{graphicx} 

\definecolor{Light}{gray}{0.85}

\def\abs#1{\left\vert #1 \right\vert}
\def\allpoly{\mbox{$\re\langle X \rangle$}}
\def\allpolyell{\mbox{$\re^{\ell}\langle X \rangle$}}
\def\allpolyx0degn{\mbox{$P_n$}}

\def\allseries{\mbox{$\re\langle\langle X \rangle\rangle$}}

\def\allseriesclousureGCell{\mbox{$\overbar{\re^\ell_{GC}\langle\langle X \rangle\rangle}$}}

\def\allseriesell{\mbox{$\re^{\ell} \langle\langle X \rangle\rangle$}}

\def\allseriesellGevreys{\mbox{$\re^\ell_{\gamma}\langle\langle X \rangle\rangle$}}

\def\allseriesLC{\mbox{$\re_{LC}\langle\langle X \rangle\rangle$}}

\def\allseriesellLC{\mbox{$\re^{\ell}_{LC}\langle\langle X \rangle\rangle$}}
\def\allseriesellGC{\mbox{$\re^{\ell}_{GC}\langle\langle X \rangle\rangle$}}

\def\allseriesX1{\mbox{$\re [[ X_1 ]]$}}
\def\Bpm{B_{\mathfrak{p}}^m}

\def\Bqell{B_{\mathfrak{q}}^\ell}

\def\bull{\rule{0.08in}{0.08in}} 

\newcommand{\comment}[1]{} 

\def\eqref#1{(\ref{#1})} 

\def\Lpm{L_{\mathfrak{p}}^m}

\def\Lpme{L^m_{\mathfrak{p},e}}
\def\LpGerm{L_{\mathfrak{p}}^m}

\def\LqGermell{L_{\mathfrak{q}}^\ell}

\def\nat{{\mathbb N}} 
\def\norm#1{\Vert#1\Vert}

\newcommand{\overbar}[1]{\mkern 1mu\overline{\mkern-1mu#1\mkern-3mu}\mkern 3mu}

\def\re{{\mathbb R}} 

\def\shuffle{{\scriptscriptstyle \;\sqcup \hspace*{-0.05cm}\sqcup\;}}

\def\begce{\begin{center}}
\def\endce{\end{center}}
\def\begar{\begin{array}}
\def\endar{\end{array}}
\def\begeq{\begin{equation}}
\def\endeq{\end{equation}}
\def\begdi{\begin{displaymath}}
\def\enddi{\end{displaymath}}
\def\begdis{\begin{eqnarray*}}
\def\enddis{\end{eqnarray*}}
\def\begeqa{\begin{eqnarray}}
\def\endeqa{\end{eqnarray}}
\def\begdes{\begin{description}}
\def\enddes{\end{description}}
\def\begit{\begin{itemize}}
\def\endit{\end{itemize}}
\def\begen{\begin{enumerate}}
\def\enden{\end{enumerate}}
\def\beglar{\left[\begin{array}}
\def\endrar{\end{array}\right]}
\def\begle{\begin{lem}}
\def\endle{\end{lem}}
\def\begde{\begin{defn}}
\def\endde{\end{defn}}
\def\begth{\begin{thm}}
\def\endth{\end{thm}}
\def\begco{\begin{cor}}
\def\endco{\end{cor}}
\def\begprop{\begin{prop}}
\def\endprop{\end{prop}}
\def\begex{\begin{example}}
\def\endex{\end{example}}
\def\begexer{\begin{exercise}}
\def\endexer{\end{exercise}}

\def\begres{\noindent{\bf Remarks}:\begin{enumerate}}
\def\endres{\end{enumerate} \par}
\def\begpr{\noindent{\em Proof:}$\;\;$}
\def\endpr{\hfill\bull \vspace*{0.05in}}
\def\begtab{\begin{tabular}}
\def\endtab{\end{tabular}}
\def\rref#1{(\ref{#1})}

 \newcommand{\ellInfty}[3]{\ell_{\infty,#1}(#2,#3)}
 \newcommand{\ellInftyNorm}[2]{\norm{#1}_{\ell_{\infty,#2}}}
 
 \newcommand{\R}{\mathbb{R}}
 \newcommand{\Frechet}{Fr\'echet}

\DeclareMathOperator{\supp}{supp}
\DeclareMathOperator{\dist}{dist}
\newcommand{\ord}{\mathrm{ord}}
\newcommand{\charseries}{\mathrm{char}}

\allowdisplaybreaks[4]
\newtheorem{example}[subsection]{Example}
\newtheorem{prop}[subsection]{Proposition}
\newtheorem{thm}[subsection]{Theorem}
\newtheorem{lem}[subsection]{Lemma}

\begin{document}

\title{Continuity of Chen-Fliess Series for Applications in
System Identification and Machine Learning}

\author{Rafael Dahmen\footnote{Karlsruhe Institute of Technology, 76131 Karlsruhe, Germany, (e-mail: rafael.dahmen@kit.edu)}, W.\ Steven Gray \footnote{Old Dominion University, Norfolk, Virginia 23529 USA, (e-mail: sgray@odu.edu)}, Alexander Schmeding\footnote{Department of Mathematics, Universitet i Bergen, All\'{e}gate 41, 5020 Bergen, Norway,
(e-mail:alexander.schmeding@uib.no)}}
\date{}
\maketitle

\begin{abstract} 
Model continuity plays an important role in applications like system identification, adaptive control, and machine learning.
This paper provides sufficient conditions under which input-output systems represented by
locally convergent Chen-Fliess series are \emph{jointly} continuous
with respect to their generating series and as operators mapping a ball in an $L_{\mathfrak p}$-space to a ball in an $L_{\mathfrak q}$-space,
where $\mathfrak p$ and $\mathfrak q$ are conjugate exponents.
The starting point is to introduce a class of topological vector
spaces known as \emph{Silva spaces} to frame the problem and then to employ the concept of a
\emph{direct limit} to describe convergence. The proof of the main continuity result combines elements of proofs for other
forms of continuity appearing
in the literature to produce the desired conclusion.
\end{abstract}
\smallskip 

\tableofcontents

\section{Introduction}

In applications involving system identification, adaptive control, and machine learning,
a stream of input-output data is continually processed over time to produce a sequence of parameter/weight estimates
so that an
assumed model's behavior matches that of the data source.
In the context of control, for example, this usually means that the dynamics of the model should asymptotically approach those of the plant. This can fail to happen
when the model is incompatible with the plant or the data stream contains insufficient information.
A more subtle mode of failure is one where the model's dynamics do not depend \emph{continuously}
on the parameters. In which case, it is possible for the sequence of parameter estimates to converge to some limit, while the corresponding sequence of approximations of the model's dynamics fail to converge in any sense.

The earliest work on the continuity of input-output systems was that of Hazewinkel \cite{Hazewinkel_80}. The focus
there was on one parameter families of linear time-invariant systems and certain \emph{degeneration} phenomena.
Continuity of the same class of systems was later address from the behaviorial point of view in
\cite{Nieuwenhuis-Willems_88,Nieuwenhuis-Willems_92}. Continuity of one parameter families of input-output systems
with Chen-Fliess series representations \cite{Fliess_81} was first
characterized by \cite{Wang_90}. In this same work it was also shown that under certain growth conditions on the generating
series such system are continuous as maps from
$L_1[0,T]$ into $C[0,T]$ with the $L_\infty$-norm for $T>0$ sufficiently small. More stringent growth conditions can even render
an output function which is well defined and continuous on $[0,\infty)$ \cite{Gray-Wang_02}. Various improvements and generalizations
of these result have
appeared in \cite{Duffaut_Espinosa_09,Winter_Arboleda_19}.
In parallel with this development,
continuity properties regarding control affine nonlinear
state space models have appeared in \cite{Azhmyakov_09}. The primary aim there was to characterize the continuity of flows with respect to the
input and initial condition. Continuity with respect to the vector fields of the realization was not considered. As the coefficients of the
corresponding Chen-Fliess depend explicitly on these vector fields and the initial condition, that analysis will not directly apply here.

The main objective of this paper is provide sufficient conditions under which input-output systems represented by
locally convergent Chen-Fliess series are \emph{jointly} continuous
with respect to their generating series and as operators mapping a ball in an $L_{\mathfrak p}$-space to a ball in an $L_{\mathfrak q}$-space,
where $\mathfrak p$ and $\mathfrak q$ are conjugate exponents.
Of course, continuity and convergence are ultimately topological concepts, so this
phenomenon can only be understood precisely in a topological framework.
The starting point is to introduce a class of topological vector
spaces known as \emph{Silva spaces} to frame the problem and then to employ the concept of a
\emph{direct limit} to describe convergence.
The proof of the main continuity result combines elements of proofs for weaker
forms of continuity appearing
in \cite{Wang_90}, \cite{Gray-Wang_02}, and \cite{Duffaut_Espinosa_09} to produce the desired conclusion.

The paper is organized as follows. The next section gives a brief summary of the Chen-Fliess series mainly to establish the notation.
The subsequent section describes the topological concepts used throughout the paper. The main continuity results appear in Section~\ref{sec:continuity-results}
along with some examples to illustrate their application.
The final section summarizes the paper's main conclusions.
\smallskip

\textbf{Acknowledgment:} The second author was supported by the National Science
Foundation under grant CMMI-1839378.

\section{Chen-Fliess Series}

An \emph{alphabet} $X=\{ x_0,x_1,$ $\ldots,x_m\}$ is any nonempty and finite set
of noncommuting symbols referred to as \emph{letters}. A \emph{word} $\eta=x_{i_1}\cdots x_{i_k}$ is a finite sequence of letters from $X$.
The number of letters in a word $\eta$, written as $\abs{\eta}$, is called its \emph{length}.
The empty word, $\emptyset$, is taken to have length zero.
The collection of all words having length $k$ is denoted by
$X^k$. Define $X^\ast=\bigcup_{k\geq 0} X^k$ and $X^{\leq J}=\bigcup_{k = 0}^J X^k$.
The former is a monoid under the concatenation product.
Any mapping $c:X^\ast\rightarrow
\re^\ell$ is called a \emph{formal power series}.
Often $c$ is
written as the formal sum $c=\sum_{\eta\in X^\ast}(c,\eta)\eta$,
where the \emph{coefficient} $(c,\eta)$ is the image of
$\eta\in X^\ast$ under $c$.
The \emph{support} of $c$, $\supp(c)$, is the set of all words having nonzero coefficients.
The set of all noncommutative formal power series over the alphabet $X$ is
denoted by $\allseriesell$. The subset of series with finite support, i.e., polynomials,
is represented by $\allpolyell$.
Each set is an associative $\re$-algebra under the catenation product and an associative and commutative $\re$-algebra under the \emph{shuffle product}, that is, the bilinear product uniquely specified by the shuffle product of two words
\begdi
(x_i\eta)\shuffle(x_j\xi)=x_i(\eta\shuffle(x_j\xi))+x_j((x_i\eta)\shuffle \xi),
\enddi
where $x_i,x_j\in X$, $\eta,\xi\in X^\ast$ and with $\eta\shuffle\emptyset=\emptyset\shuffle\eta=\eta$ \cite{Fliess_81}.

Given any $c\in\allseriesell$ one can associate a causal
$m$-input, $\ell$-output operator, $F_c$, in the following manner.
Let $\mathfrak{p}\ge 1$ and $t_0 < t_1$ be given. For a Lebesgue measurable
function $u: [t_0,t_1] \rightarrow\re^m$, define
$\norm{u}_{\mathfrak{p}}=\max\{\norm{u_i}_{\mathfrak{p}}: \ 1\le
i\le m\}$, where $\norm{u_i}_{\mathfrak{p}}$ is the usual
$L_{\mathfrak{p}}$-norm for a measurable real-valued function,
$u_i$, defined on $[t_0,t_1]$.  Let $L^m_{\mathfrak{p}}[t_0,t_1]$
denote the set of all measurable functions defined on $[t_0,t_1]$
having a finite $\norm{\cdot}_{\mathfrak{p}}$ norm and
$B_{\mathfrak{p}}^m(R_u)[t_0,t_1]:=\{u\in
L_{\mathfrak{p}}^m[t_0,t_1]:\norm{u}_{\mathfrak{p}}\leq R_u\}$.
Assume $C[t_0,t_1]$
is the subset of continuous functions in $L_{1}^m[t_0,t_1]$. Define
inductively for each $\eta\in X^{\ast}$ the map $E_\eta:
L_1^m[t_0, t_1]\rightarrow C[t_0, t_1]$ by setting
$E_\emptyset[u]=1$ and letting
\[E_{x_i\bar{\eta}}[u](t,t_0) =
\int_{t_0}^tu_{i}(\tau)E_{\bar{\eta}}[u](\tau,t_0)\,d\tau, \] where
$x_i\in X$, $\bar{\eta}\in X^{\ast}$, and $u_0=1$. The
\emph{Chen-Fliess series} corresponding to $c$ is
\begeq
y(t)=F_c[u](t) =
\sum_{\eta\in X^{\ast}} (c,\eta)\,E_\eta[u](t,t_0) \label{eq:Fliess-operator-defined}
\endeq
\cite{Fliess_81,Fliess_83}.
It can be shown that if there exists real numbers $K,M\geq 0$ such that
\begeq \label{eq:LC-condition}
\abs{(c,\eta)}\leq KM^{\abs{\eta}}\abs{\eta}!,\;\;\forall \eta\in X^\ast
\endeq
($\abs{z}:=\max_i \abs{z_i}$ when $z\in\re^\ell$)
then the series defining $F_c$ converges absolutely and uniformly for sufficient
small $R,T >0$ and
constitutes a well defined mapping from
$B_{\mathfrak p}^m(R)[t_0,$ $t_0+T]$ into $B_{\mathfrak
q}^{\ell}(S)[t_0, \, t_0+T]$, where the numbers $\mathfrak{p},\mathfrak{q}\in[1,+\infty]$ are
conjugate exponents, i.e., $1/\mathfrak{p}+1/\mathfrak{q}=1$ \cite{Gray-Wang_02}.
Any such mapping is called a \emph{locally convergent Fliess operator}.

A more refined convergence analysis of Chen-Fliess series appears in \cite{Winter_Arboleda_19}
utilizing the notion of \emph{Gevrey order}. A series $c \in \allseriesell$ is said to have Gevrey order $s \in [0,\infty)$ if there
exists constants $K,M>0$ such that
\begeq
\abs{(c,\eta)}\leq K M^{\abs{\eta}} (\abs{\eta}!)^s,\;\; \forall \eta\in
X^\ast. \label{eq:Gevrey-growth-condition}
\endeq
Clearly, if $c$ has Gevrey order $s$ then it is also has Gevrey order $s^ \prime$, where $s^ \prime > s$.
Define for a given $c$ the real number $\gamma_c=\min\{s \in [0,\infty): s  \text{ satisfies~\eqref{eq:Gevrey-growth-condition}}\}$ and
the set of all generating series with minimum Gevrey order $\gamma$ as $\allseriesellGevreys$.
In this context, the set of all generating series for locally convergent Fliess operators as described above
is
\begdi
\allseriesellLC:= \bigcup_{0 \leq \gamma \leq 1} \allseriesellGevreys,
\enddi
while a smaller set of series (note the upper bound on $\gamma$)
\begdi
\allseriesellGC:= \bigcup_{0 \leq \gamma < 1} \allseriesellGevreys,
\enddi
can be shown to yield a type of \emph{global convergence}
on the extended space
$\Lpme (t_0)$ into $C[t_0, \infty)$,  where
\begdi
\Lpme(t_0):=\bigcup_{T>0}\Lpm[t_0,t_0+T]
\enddi
\cite{Winter_Arboleda-etal_2015}. Interestingly, this latter set of generating series does
\emph{not} constitute all of those that provide a globally defined Fliess operator as shown by
example in \cite{Winter_Arboleda_19}.

Finally, a Fliess operator $F_c$ defined on $B_{\mathfrak p}^m(R)[t_0,t_0+T]$
is said to be \emph{realizable} when there exists
a state space model
\begin{subequations} \label{eq:general-MIMO-control-affine-system}
\begin{align}
\dot{z}(t)&= g_0(z(t))+\sum_{i=1}^m g_i(z(t))\,u_i(t),\;\;z(t_0)=z_0  \label{eq:state}\\
y_j(t)&=h_j(z(t)),\;\;j=1,2,\ldots,\ell, \label{eq:output}
\end{align}
\end{subequations}
where each $g_i$ is an analytic vector field expressed in local
coordinates on some neighborhood ${\cal W}$ of $z_0$,
and each output function $h_j$ is an analytic function on ${\cal W}$ such that
(\ref{eq:state}) has a well defined solution $z(t)$, $t\in[t_0,t_0+T]$ for
any given input $u\in B_{\mathfrak
p}^m(R)[t_0,t_0+T]$, and $y_j(t)=F_{c_j}[u](t)=h_j(z(t))$, $t\in[t_0,t_0+T]$, $j=1,2,\ldots,\ell$.
It can be shown that for any word $\eta=x_{i_k}\cdots x_{i_1}\in X^\ast$
\begeq
(c_j,\eta)=L_{g_{\eta}}h_j(z_0):=L_{g_{i_1}}\cdots L_{g_{i_k}}h_j(z_0), \label{eq:c-equals-Lgh}
\endeq
where $L_{g_i}h_j$ is the \emph{Lie derivative} of $h_j$ with respect to $g_i$.

\section{Topological subspaces of $\allseriesell$}

Suppose a sequence of generating series $\{c_j\}_{j\geq 1}$ is produced in real-time
by processing a stream of input-output data in some manner. The
corresponding sequence of Chen-Fliess series is taken to be $\{F_{c_j}\}_{j\geq 1}$.
If the estimation or learning algorithm producing these generating series ensures that $c_j\rightarrow c$ in some sense,
then it is desirable that $F_{c_j}\rightarrow F_c$ in some
fashion as well. Perhaps the most obvious way in which one series can approach another is in the ultrametric sense.
Specifically, for any fixed real number $\sigma$ such that $0<\sigma<1$,
consider the mapping
\begin{align*} \label{eq:dist}
\dist&:\allseries \times \allseries \rightarrow \re,\\
      &\ (c,d)\mapsto \sigma^{\ord(c-d)},
\end{align*}
where $\ord(c)$ is the length of the shortest word in the support of $c$ ($\ord(0):=\infty$).
The $\re$-vector space $\allseriesell$ with mapping $\dist$ is known to be a complete ultrametric space \cite{Berstel-Reutenauer_88}.
If each series $c_j\in\allseriesLC$, the following simple example illustrates that in the limit there is not always a well defined operator
to which a given sequence of Fliess operators is converging.

\begex \label{exa:noconv}
\textrm{
Let $X=\{x_1\}$ and
consider the sequence of polynomials}
\begdi
c_j=x_1+(2!)^2\;x_1^2+(3!)^2\;x_1^3+\cdots + (j!)^2\; x_1^j,\;\; j\geq 1.
\enddi
\textrm{Clearly, each polynomial $c_j$ is locally convergent. Thus,
each Fliess operator $F_{c_j}$ is well defined
on some ball of input functions in $L^m_{\mathfrak{p}}[t_0,t_1]$.
Furthermore,
the sequence $(c_j)_j$ converges to $c=\sum_{k\geq 1} (k!)^2x_1^k$ in the ultrametric topology.
But the limiting Chen-Fliess series $F_c$ is not
well defined in any obvious sense.
}
\endex

This example motivates the following fundamental problem: On what topological subspaces of $\allseriesell$ does $c_j\rightarrow c$ imply that $F_{c_j}\rightarrow F_{c}$ in some sense with the limit point $F_c$ being a well defined operator?
The following subsections lay the foundation for addressing this problem by presenting what subspaces are available for consideration.

\subsection{Fixed $M>0$ (Banach Spaces)}

As a first step, consider the following interpretation of condition (\ref{eq:LC-condition}). Fix $M>0$ and define
\[
       \ellInftyNorm{c}{M} := \sup\left\{   \frac{\abs{(c,\eta)}}{M^{\abs{\eta}}\abs{\eta}!} : \eta\in X^\ast    \right\} \in \left[0,\infty\right]
\]
for each $c\in\allseriesell$.
The set of all $c$ with $\ellInftyNorm{c}{M}<\infty$ is denoted by $\ellInfty{M}{X^*}{\R^\ell}$.
It is straightforward to check that $\ellInfty{M}{X^*}{\R^\ell}$
is a vector subspace of $\allseriesell$. The function $\ellInftyNorm{\cdot}{M}$ is a norm on $\ellInfty{M}{X^*}{\R^\ell}$.
The following assignment is an isometry of normed spaces:
\[
\ellInfty{M}{X^*}{\R^\ell}\longrightarrow \ell_\infty(X^*,\R^\ell)
: c \mapsto \frac{c}{M^{\abs{\eta}}\abs{\eta}!},
\]
where
$\ell_\infty(X^*,\R^\ell):= \left\{c \colon X^* \rightarrow \R^\ell : \sup_{\eta }\abs{(c,\eta)} <\infty\right\}$ is the Banach space of all \emph{bounded} functions from $X^*$ to $\R^\ell$.
This shows that for each fixed $M>0$ the space $(\ellInfty{M}{X^*}{\R^\ell},\ellInftyNorm{\cdot}{M})$
is a Banach space.
A series $c\in\allseriesell$ belongs to $\ellInfty{M}{X^*}{\R^\ell}$ if and only if the bound (\ref{eq:LC-condition}) holds for some $K\geq0$ and the fixed number $M$. In fact, the norm
$\ellInftyNorm{c}{M}$ is the smallest number $K\geq0$ such that (\ref{eq:LC-condition}) is satisfied.

As $\ellInfty{M}{X^*}{\R^\ell}$ is a Banach space, and, in particular, a metric space,
the topology of $\ellInfty{M}{X^*}{\R^\ell}$ can be recovered from convergent sequences, where a sequence $(c_j)_j$ in $\ellInfty{M}{X^*}{\R^\ell}$ converges to $c\in\ellInfty{M}{X^*}{\R^\ell}$ if and only if
\[
   \lim_{j\to\infty}\ellInftyNorm{c_j - c}{M} =0.
\]
Given that $\ellInfty{M}{X^*}{\R^\ell}$ is an infinite dimensional Banach space, the Bolzano-Weierstrass theorem fails to hold, i.e., not every $\ellInftyNorm{\cdot}{M}$-bounded
sequence has a $\ellInftyNorm{\cdot}{M}$-convergent subsequence, see \cite[Satz I.2.7]{Werner}.
Furthermore, the space is not separable, i.e., there is no countable dense subset.
Given $M_1$ and $M_2$ such that $M_1 \leq M_2$, it is clear that
$ \ellInftyNorm{\cdot}{M_1} \geq \ellInftyNorm{\cdot}{M_2}$,
and thus the inclusion (as vector spaces)
\[
   \ellInfty{M_1}{X^*}{\R^\ell} \subseteq \ellInfty{M_2}{X^*}{\R^\ell}
\]
holds.
This inclusion is not a topological embedding as the topology induced by $\ellInfty{M_2}{X^*}{\R^\ell}$ is coarser than the one induced by $\ellInfty{M_1}{X^*}{\R^\ell}$.
It turns out for $M_1 < M_2$ that the inclusion map
\begeq \label{eq:inclusion maps}
   \ellInfty{M_1}{X^*}{\R^\ell} \to \ellInfty{M_2}{X^*}{\R^\ell}
\endeq
is a \emph{compact operator} (see
\cite[Lemma B.6]{DS18}), i.e., it maps bounded sets to relatively compact sets. In particular, this shows for $M_1<M_2$ that
every sequence which is bounded in the $ \ellInftyNorm{\cdot}{M_1}$-norm has a subsequence which converges in the coarser $ \ellInftyNorm{\cdot}{M_2}$-topology.

\subsection{The projective limit $M\to0$ (\Frechet{}--Schwartz spaces)}

Consider next those $c\in\allseriesell$ for which $\ellInftyNorm{c}{M}$ is finite for all $M>0$.
This means that for each $M>0$ there is a $K=\ellInftyNorm{c}{M}\geq0$ satisfying \rref{eq:LC-condition}.
Algebraically, this corresponds to the \emph{intersection} of all vector spaces $\ellInfty{M}{X^*}{\R^\ell}$, namely,
\[
   \ellInfty{\leftarrow}{X^*}{\R^\ell} := \bigcap_{M>0} \ellInfty{M}{X^*}{\R^\ell}.
\]
On spaces like these, there is a natural topology which turns this space into a locally convex topological vector space.
In the functional analysis literature, this object is called the \emph{projective limit} (or \emph{inverse limit} or \emph{categorical limit}) of the system $\left(       \ellInfty{M}{X^*}{\R^\ell} \right)_{M>0}$ and denoted also by
\begin{align*}
   \ellInfty{\leftarrow}{X^*}{\R^\ell}:=& \lim_{\substack{\longleftarrow \\M\to0 }}   \ellInfty{M}{X^*}{\R^\ell}   = \bigcap_{M>0} \ellInfty{M}{X^*}{\R^\ell}.
\end{align*}
For a given $c\in\allseriesell$, one can check whether it belongs to this space in the following way:
\begin{equation}
 c\in\ellInfty{\leftarrow}{X^*}{\R^\ell} \\\Longleftrightarrow \ellInftyNorm{c}{M}<\infty,\; \forall M>0.
\label{eq:proj_condition}
\end{equation}
The sequence $M_k = 1/k$, $k\in \mathbb{N}$ is cofinal, hence it suffices to check \rref{eq:proj_condition} only for $M$ of the form $M_k$.

Now $\ellInfty{\leftarrow}{X^*}{\R^\ell}$ is the projective limit of countably many
Banach spaces. Thus it becomes a \emph{\Frechet{} space}, i.e.\ a complete metrisable space. \Frechet{} spaces share many nice properties with Banach spaces, e.g.,
their topology is determined by sequences,
where a sequence $(c_j)_j$ in $\ellInfty{\leftarrow}{X^*}{\R^\ell}$ converges to $c\in\ellInfty{\leftarrow}{X^*}{\R^\ell}$ if and only if
\[
  \lim_{j\to\infty}\ellInftyNorm{c_j - c}{M} =0,\;\forall M>0.
\]
(Again it suffices to check this only for all $M=1/k$, $k\in\mathbb N$.)
Since the inclusion maps are all compact operators, $\ellInfty{\leftarrow}{X^*}{\R^\ell}$ is even a \Frechet{}--Schwartz space, \cite[Definition 8.5.2]{BaP87}. Hence, it behaves
much nicer than the Banach spaces from which it was build.
In particular, the space $\ellInfty{\leftarrow}{X^*}{\R^\ell}$ satisfies a version of the Bolzano-Weierstrass theorem, namely, every $\ell_{\infty,\leftarrow}$-bounded sequence has a $\ell_{\infty,\leftarrow}$-convergent subsequence. Here, a sequence $(c_j)_j$ is called $\ell_{\infty,\leftarrow}$-bounded if for \emph{all} $M>0$
it holds that
$\sup_{j}\ellInftyNorm{c_j}{M} < \infty$.
This follows from \cite[Proposition 8.5.9]{BaP87}, which furthermore implies that $\ellInfty{\leftarrow}{X^*}{\R^\ell}$ is separable, i.e., there is countable dense subset.
In \cite{Winter_Arboleda_19}[Theorem 3.4.5] it is shown that
\begdi
\ellInfty{\leftarrow}{X^*}{\R^\ell}=\allseriesclousureGCell,
\enddi
where the closure on the right is taken with respect to the $\ell_{\infty,\leftarrow}$-topology (called the \emph{semi-norm topology} in loc.cit.).
In other words, there are some generating series with minimum Gevrey order $\gamma=1$ that
yield globally defined Fliess operators.

\subsection{The direct limit $M\to\infty$ (Silva Spaces)}
Consider next a series $c\in\allseriesell$ where there exists at least one number $M>0$ such that $K=\ellInftyNorm{c}{M}\geq0$ satisfies (\ref{eq:LC-condition}).
Algebraically, this case corresponds to the \emph{union} of all vector spaces $\ellInfty{M}{X^*}{\R^\ell}$, that is, \[
  \ellInfty{\rightarrow}{X^*}{\R^\ell} := \bigcup_{M>0} \ellInfty{M}{X^*}{\R^\ell}.
\]
As with the intersection, there is also a natural topology turning this space into a locally convex topological vector space. This object is
called the \emph{direct limit} (or \emph{inductive limit} or \emph{categorical colimit}) of the system
$\left(\ellInfty{M}{X^*}{\R^\ell} \right)_{M>0}$ and denoted by
\begin{align*}
  \ellInfty{\rightarrow}{X^*}{\R^\ell} =& \lim_{\substack{\longrightarrow \\M\to\infty }}   \ellInfty{M}{X^*}{\R^\ell} \\ =&  \bigcup_{M>0} \ellInfty{M}{X^*}{\R^\ell}.
\end{align*}
This construction can also be found in \cite{BS16,DS18}.
For a given $c\in\allseriesell$, one can check whether it belongs to this space in the following way:
\begin{equation}
  c\in\ellInfty{\leftarrow}{X^*}{\R^\ell} \Longleftrightarrow \exists M>0 \text{ such that}\; \ellInftyNorm{c}{M}<\infty.   \label{eq:direct_condition}
\end{equation}
Since the sequence $M_k = k$, $k\in \mathbb N$ is cofinal, one can always find an $M\in\mathbb N$ for which $\ellInftyNorm{c}{M}<\infty$.
Thus, one could equivalently work only with $M\in\mathbb N$.
In general, direct limits are more difficult to work with than projective limits. Fortunately, this particular direct limit is a countable direct limit of Banach spaces with compact operators are inclusion maps. Direct limit spaces like these are called \emph{Silva space}.

Although Silva spaces are not metrizable, they are always \emph{sequential} \cite[Proposition 6]{Yos57}. This means that as in the Banach space case, the topology is determined by sequences, i.e.\ sets are closed if and only if they are sequentially closed. A sequence $(c_j)_j$ in $\ellInfty{\rightarrow}{X^*}{\R^\ell}$ converges to  $c\in\ellInfty{\rightarrow}{X^*}{\R^\ell}$ if and only if
\[
  \lim_{j\to\infty}\ellInftyNorm{c_j - c}{M} =0\text{ for \emph{one fixed} }M>0.
\]
In other words, a sequence in a Silva space converges if there exists one fixed $M>0$ for the whole sequence such that $(c_j)_j$ converges in the Banach space $\ellInfty{M}{X^*}{\R^\ell}$ \cite[Theorem 1]{Yos57}.
In particular, note that for a sequence to converge, all terms must lie in one of the spaces $\ellInfty{M}{X^*}{\R^\ell}$, i.e., one $M>0$ has to work for the whole sequence.
The sequence in Example \ref{exa:noconv} fails to converge in the Silva topology since there is no $M$ for which the sequence is Cauchy.

Using again the compactness of the inclusion maps, it follows that a sequence which is bounded in one $ \ellInftyNorm{\cdot}{M_1}$-norm (for a given $M_1>0$) has a subsequence which converges in the coarser $ \ellInftyNorm{\cdot}{M_2}$-topology for all $M_2>M_1$. As earlier, there is a version of the Bolzano-Weierstrass theorem, namely, every $\ell_{\infty,\rightarrow}$-bounded sequence has a $\ell_{\infty,\rightarrow}$-convergent subsequence.
In this case, a sequence $(c_j)_j$ is called $\ell_{\infty,\leftarrow}$-bounded if there is at least one $M>0$ with $\sup_{j}\ellInftyNorm{c_j}{M} < \infty$.
Therefore, a Silva space has better topological properties than the Banach spaces from which it is constructed.
Furthermore, every Silva space is separable.
Finally, it is shown in \cite[Theorem 3.2.7]{Winter_Arboleda_19} that
\begdi
\ellInfty{\rightarrow}{X^*}{\R^\ell}=\allseriesellLC.
\enddi

\begex \label{ex:Bananch-convergence}
Let $X=\{x_1\}$. The sequence $c_j := j!x_1^j$, $j\in \mathbb{N}$ has norm
$\lVert c_j \rVert_{\ell_{\infty,1}} = 1$. Therefore, $(c_j)_j$ does not converge to zero in $\ell_{\infty,1} (X^\ast, \R)$.
However, since $\lVert c_j\rVert_{\ell_{\infty,2}} = 1/2^{j}$, it follows that $c_j \xrightarrow{j\rightarrow \infty} 0$ in $\ell_{\infty,2} (X^\ast,\R)$ and also in the Silva topology.
\endex

\begex \label{ex:Silva-convergence}
Define for $n, j\in \mathbb{N}$ the sequence $d_{n,j}:=n^{(5j-2)/2j}C_n$, where $C_n := (2n)!/((n+1)!n!)$ is the $n$th Catalan number.\footnote{Sequence A000108 in OEIS.} Recall that the asymptotic growth of the Catalan numbers is $C_n \sim 4^n/(n^{3/2}\sqrt{\pi})$. Thus, for $d_j := \sum_{n=1}^\infty d_{n,j} x_1^{n}$, it is clear that $d_1 \in \ell_{\infty,4} (X^*, \R)$, but $d_j \not \in \ell_{\infty,4} (X^*, \R)$ for $j>1$. However, since $\lVert d_j \rVert_{\ell_{\infty,5}} < \infty$, it does hold that $(d_j)_j \subseteq \ell_{\infty,\rightarrow} (X^*,\R)$. Furthermore, it is
easily checked that $\lim_{j\rightarrow \infty}\lVert d_j -d \rVert_{\ell_{\infty,5}}=0$, where $d= \sum_{n=1}^\infty n^{5/2}C_{n}x^{n}_1$.
Thus, $d_j$, $j\in\mathbb N$ converges to zero in the Silva topology.
\endex

In the continuity theorems presented in the next section, every sequences $(c_j)_j$ will be assumed a priori to be entirely contained in some Banach space $\ell_{\infty,M} (X^*, \R)$, $M>0$, thus avoiding the phenomenon shown in the previous example. Therefore, only a Banach topology is really needed. However, for applications
such as the interconnection of Chen-Fliess operators, the Silva topology is more applicable in the corresponding continuity analysis.
For example, one can define for generating series $c,d \in \ell_{\infty,M} (X^*,\R^\ell)$ a product $c \circ d$ such that the composition
satisfies $F_c\circ F_d=F_{c\circ d}$. It is well known that in general $ c \circ d$ will not be contained in $\ell_{\infty,M} (X^*,\R^\ell)$.
However, there exists a $K(M,N) < \infty$ for all $M,N \in [0,\infty[$ such that
$$\circ\colon \ell_{\infty, M} (X^*,\R^\ell) \times \ell_{\infty, N} (X^*,\R^\ell)  \rightarrow \ell_{\infty, K(M,N)} (X^*,\R^\ell) $$
is well defined \cite{Thitsa-Gray_SIAM12}. Using these estimates, the interconnection of Chen-Fliess series induces a continuous product on the Silva space $\ell_{\infty, \rightarrow} (X^*,\R^\ell)$. Hence, the Silva topology is the natural topology for describing the continuity of such interconnections.

\section{Main Continuity Theorems}
\label{sec:continuity-results}

The continuity problem for a Chen-Fliess series $F_c[u]$ is approached incrementally. It is first assumed that
the input $u$ is fixed and the generating series $c$ is variable (series to output continuity). Then the case where $c$ is fixed and $u$ is variable is presented (input-output operator continuity). Finally, the two cases are combined.
For notational convenience, define the space of $L_{\mathfrak{p}}$-germs
$\LpGerm (t_0) := \{[u] \mid u \in L_{\mathfrak{p}}^m[t_0,t_1] \text{ for some } t_1> t_0\}$,
where the class $[u]$ contains all functions equal to $u$ in some neighborhood of $t_0$.
Note that this space can \emph{not} be endowed with any useful topology making the inclusion
$L_{\mathfrak{p}}^m[t_0,t_1]$ continuous (as this would automatically be non-Hausdorff).

\begth[series to output continuity]\label{thm:outputcont}
The map
$$\LpGerm(0) \times \allseriesellLC \rightarrow \LqGermell(0),\quad (u,c)\mapsto y=F_c[u]$$
is well defined. Moreover, for every $M>0$ and fixed $u\in \LpGerm(0)$, there exists a $T>0$ such that
$$\ell_{\infty,M} (X^*,\R^\ell) \rightarrow L^\ell_{\mathfrak{q}}[0,T],\quad c \mapsto y=F_c[u]$$
is continuous.
\endth

The following two lemmas are needed for the proof \cite{Duffaut_Espinosa_09}.

\begle \label{le:DE1}
Let $X=\{x_0,x_1,\ldots,x_m\}$. For any $k\in\nat_0$, the characteristic polynomial $\charseries (X^k)$ of $X^k$, i.e.~the sum of all words of length $k$, satisfies
\begeq
\charseries(X^k)
= \sum_{\genfrac{}{}{0pt}{2}{r_0,r_1,\ldots,r_m\geq 0}{ r_0+r_1+\cdots+r_m=k}}
x_0^{r_0}\shuffle x_1^{r_1}\shuffle \cdots \shuffle x_m^{r_m},
\label{eq:Luis-charXk-identity}
\endeq
where $\shuffle$ denotes the shuffle product.
\endle

\begle\label{le:DE2}
Let $X=\{x_0,x_1,\ldots,x_m\}$. For any $u\in L_1^m[0,T]$ and $\eta\in X^\ast$
\begdi
\abs{E_\eta[u](t)}\leq E_\eta[\bar{u}](t),\;\; 0\leq t\leq T,
\enddi
where $\bar{u}\in L_1^m[0,T]$ has components $\bar{u}_j:=\abs{u_j}$,
$j=1,2,\ldots,m$.
Furthermore, for any integers $r_j\geq 0$ it follows that
\begdi
\abs{E_{x_0^{r_0}\shuffle x_1^{r_1}\shuffle\cdots\shuffle x_m^{r_m}}[u](t)}
\le \prod_{j=0}^m \frac{U_j^{r_j}(t)}{r_j!},
\;\;0\leq t\leq T,
\enddi
where
$U_j(t) := \int_0^t\abs{u_j(\tau)}\,d\tau$.\footnote{For notational
convenience, occasionally $F_p$ will denoted by $E_p$
when $p\in\allpoly$.}
In particular, if on $[0,T]$ it is assumed that $\max\{\norm{u}_1,T\}\le R$
then
\begdi
\abs{E_{x_0^{r_0}\shuffle x_1^{r_1}\shuffle\cdots\shuffle x_m^{r_m}}[u](t)}\le
\frac{R^{k}}{r_0!\,r_1!\,\cdots r_m!},
\;\;0\leq t\leq T,
\enddi
where $k=\sum_{j}r_j$.
\endle

Now the proof of Theorem \ref{thm:outputcont}.

\begpr
If $c\in\allseriesellLC$ then there exists $K,M\geq 0$
satisfying \rref{eq:LC-condition}.
Fix $u\in L_1^m(0)$ (without loss of generality $\mathfrak{p}=1$ and $t_0=0$)
so that for some $T >0$, $u\in L_1^m[0,T]$. Define
$R = \max\{\norm{u}_1, T\}$.
Applying Lemmas~\ref{le:DE1} and \ref{le:DE2} it then follows that:
\begin{align}
\abs{y(t)}&\leq \sum_{\eta\in X^\ast} \abs{(c,\eta)E_\eta[u](t)} \leq\sum_{k=0}^\infty \sum_{\eta\in X^k} \abs{(c,\eta)} E_\eta[\bar{u}](t) \nonumber \\
&\leq\sum_{k=0}^\infty KM^kk!\sum_{\genfrac{}{}{0pt}{2}{r_0,r_1,\ldots,r_m\geq 0}{ r_0+r_1+\cdots+r_m=k}} \hspace*{-0.15in}
E_{x_0^{r_0}\shuffle x_1^{r_1}\shuffle\cdots\shuffle x_m^{r_m}}[\bar{u}](t) \nonumber\\
&\leq \sum_{k=0}^\infty KM^kk! \sum_{\genfrac{}{}{0pt}{2}{r_0,r_1,\ldots,r_m\geq 0}{ r_0+r_1+\cdots+r_m=k}} \frac{R^k}{r_0!r_1!\cdots r_m!} \nonumber\\
&=\sum_{k=0}^\infty K(MR)^k \sum_{\genfrac{}{}{0pt}{2}{r_0,r_1,\ldots,r_m\geq 0}{ r_0+r_1+\cdots+r_m=k}} \frac{k!}{r_0!r_1!\cdots r_m!} \nonumber\\
&=\sum_{k=0}^\infty K(MR(m+1))^k. \label{e-tl2}
\end{align}
Therefore, if $R < 1/(M(m+1))$, i.e., if
\begeq
\max\{\norm{u}_1, \, T\}< \frac{1}{M(m+1)}, \label{eq:T_max_inequality}
\endeq
then the series \rref{eq:Fliess-operator-defined} converges absolutely and uniformly on $[0,T]$
so that $y$ is well defined as an $L_\infty$-germ, specifically, $y\in L_\infty^m[0,T]$.
Now since the mapping $c\mapsto y=F_c[u]$ is linear, it is sufficient for proving
continuity to show that it is
bounded as a mapping from the Banach space $\ellInfty{M}{X^*}{\R^\ell}$ into
the Banach space $L_\infty^\ell[0,T]$, i.e., $\norm{y}_{\infty}/\ellInftyNorm{c}{M} < \infty$.
Observe that
$K=\ellInftyNorm{c}{M}>0$ is a valid choice, and therefore,
$\norm{y}_\infty/K\leq 1/(1-MR(m+1))<\infty$ as claimed.
\endpr

\begin{figure}
  \centering
  \includegraphics[width=8cm]{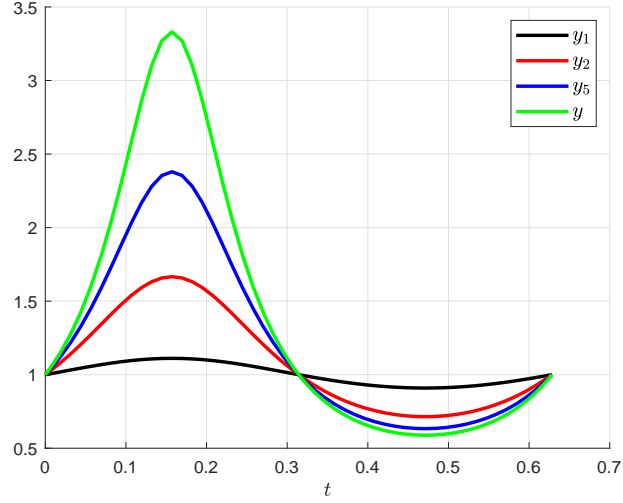}
  \caption{Convergence of $y_j$ to $y$ in Example~\ref{ex:fixed-u}
  when $M_a=1$, $M_b=7$, and $u(t)=\cos(10t)$ with $t\in[0,0.2\pi]$.}
  \label{fig:fixed-u}
\end{figure}

\begex \label{ex:fixed-u}
Let $X=\{x_1\}$ and
consider the sequence $c_j=\sum_{k\geq 0} M_j^k k!\,x_1^k$, $j\geq 1$ in $\allseriesLC$,
where
$$
M_j=M_b\theta_j+M_a(1-\theta_j)
$$
and $\theta_j=(j-1)/j$. In which case,
$
c_1=\sum_{k\geq 0} M_a^kk!\,x_1^k
$ and define
$c=\sum_{k\geq 0} M_b^k k!\,x_1^k$.
To see that $c_j \xrightarrow{j\rightarrow \infty} c$ as a direct limit, set $M=\max\{M_a,M_b\}$ and observe that
\begdi
\norm{c_j-c}_{\ell_M^\infty}=\sup_{k\geq 0}\left| \frac{M_j^k-M_b^k}{M^k}\right|=\left|\frac{M_j-M_b}{M}\right|.
\enddi
Thus, the conclusion follows directly.
Now for the given $M$ pick
$T$ and $u$ to satisfy
\rref{eq:T_max_inequality}. Therefore,
\begin{align*}
y_j(t)&=\sum_{k=0}^\infty M_j^k k! E_{x_1^k}[u](t)
=\sum_{k=0}^\infty M_j^k k! \frac{E^k[u](t)}{k!} \\
&=\frac{1}{1-M_j E_1[u](t)},
\end{align*}
and likewise
\begdi
y(t)=\frac{1}{1-M E_1[u](t)},
\enddi
are both well defined on $[0,T]$. In addition,
\begin{align*}
\|y_j-y\|_{\infty}\hspace*{-1.5pt}=\hspace*{-1.5pt}\sup_{t\in[0,T]} \left\|\frac{(M_j-M_b)E_{x_1}[u](y)}{(1-M_jE_{x_1}[u](t))(1-M_bE_{x_1}[u](t))}\right\|
\end{align*}
so that $y_j \xrightarrow{j\rightarrow \infty} y$ in the $L_\infty[0,T]$ norm sense. The specific example
where $M_a=1$, $M_b=7$, and $u(t)=\cos(10t)$ with $t\in[0,0.2\pi]$
is shown in Figure~\ref{fig:fixed-u}. Here $T=0.6283$, $\|u\|_1=0.4$ and
$1/\max(M_a,M_b)=0.1429$ (note $m+1=1$ as $X$ has only one letter), which shows that condition \rref{eq:T_max_inequality} is very
conservative in this instance.
\endex

Input-output operator continuity is addressed in the next theorem. The proof is inspired by results
appearing in \cite{Wang_90} except that certain details have to be handled differently in order
to use this result in the proof of the final continuity theorem.

\begth \label{th:Fc_operator_continuity} (input-output operator continuity)
Suppose $c\in\allseriesellLC$ and select any pair of conjugate exponents
${\mathfrak p},{\mathfrak q}\in[1,\infty]$.
If $0 < T \leq R < 1/(M(m+1))$ such that $c \in  \ellInfty{M}{X^*}{\R^\ell}$, then the operator
\begdi
F_c\colon \Bpm (R)[t_0, \, t_0+T] \rightarrow \Bqell (S)[t_0,t_0+T]
\enddi
for some $S>0$ is continuous with respect to the $L_{\mathfrak p}$ and
$L_{\mathfrak q}$ norms.
\endth

\begpr
It needs to be shown for any
$\epsilon>0$ that there exists a $\delta>0$ such that if
$v,u\in\Bpm (R)[t_0, \, t_0+T]$ satisfy
$\norm{v-u}_{\mathfrak p}<\delta$, then
$\norm{F_c[v]-F_c[u]}_{\mathfrak q}<\epsilon$.
It is first proved by induction on the length of the word $\eta\in X^\ast$ that
the mapping
\begdi
E_\eta:\Bpm (R)[t_0, \, t_0+T] \rightarrow \Bqell (S)[t_0,t_0+T]
\enddi
has the desired continuity property. The focus is on the case where
${\mathfrak p},{\mathfrak q}\in(1,\infty)$
(the remaining case is handled similarly). Without loss of generality, assume
$t_0=0$. The claim is trivial when $\eta$ is the empty word. If $\eta=x_i$, then
\begin{align*}
\norm{E_{x_i}[v]-E_{x_i}[u]}_{\mathfrak q}&=\left(\int_0^T \abs{
E_{x_i}[v](t)-E_{x_i}[u](t)
}^{\mathfrak q}dt\right)^{\frac{1}{\mathfrak q}} \\
&\leq\left(\int_0^T\left(
\int_0^T \abs{v_i(\tau)-u_i(\tau)}\,d\tau\right)^{\mathfrak q}
dt\right)^{\frac{1}{\mathfrak q}} \\
&=\int_0^T \abs{v_i(\tau)-u_i(\tau)}\,d\tau\, T^{1/{\mathfrak q}} \\
&\leq\norm{v_i-u_i}_p\,T^{2/{\mathfrak q}} \leq\norm{v-u}_p\,T^{2/{\mathfrak q}},
\end{align*}
where H\"{o}lder's inequality has been used in the second to the last step above.
Thus, if $\norm{v-u}_p<\delta_{x_i}:=\epsilon/T^{2/{\mathfrak q}}$, then clearly
\begdi
\norm{E_{x_i}[v]-E_{x_i}[u]}_{\mathfrak q}<\epsilon.
\enddi
Now suppose the claim holds for all words up to some fixed length $k\geq 0$.
Then for any $x_i\in X$ and $\eta\in X^k$ observe
\begin{align*}
\norm{E_{x_i\eta}[v]-E_{x_i\eta}[u]}_{\mathfrak q} &=\left\Vert\left(E_{x_i\eta}[v]-\int_0^\cdot u_i(\tau)E_{\eta}[v](\tau)\,d\tau\right)+\left(\int_0^\cdot u_i(\tau)E_{\eta}[v](\tau)\,d\tau-E_{x_i\eta}[u]\right)\right\Vert_{\mathfrak q} \\
&\leq\left\Vert E_{x_i\eta}[v]-\int_0^\cdot u_i(\tau)E_{\eta}[v](\tau)\,d\tau\right\Vert_{\mathfrak q}+\left\Vert\int_0^\cdot u_i(\tau)E_{\eta}[v](\tau)\,d\tau-E_{x_i\eta}[u]\right\Vert_{\mathfrak q} \\
&\leq\left(
\int_0^T\left(\int_0^T \abs{v_i(\tau)-u_i(\tau)}\abs{E_\eta[v](\tau)}\,d\tau\right)
^{\mathfrak q}dt\right)^{\frac{1}{\mathfrak q}}+ \\
&\hspace*{0.25in}\left(
\int_0^T\left(\int_0^T \abs{u_i(\tau)}\abs{E_\eta[v](\tau)-E_\eta[u](\tau)}\,d\tau\right)
^{\mathfrak q}dt\right)^{\frac{1}{\mathfrak q}} \\
&\leq \int_0^T \abs{v_i(\tau)-u_i(\tau)}\abs{E_\eta[v](\tau)}\,d\tau\, T^{1/{\mathfrak q}}+ \\
&\hspace*{0.25in} \int_0^T \abs{u_i(\tau)}\abs{E_\eta[v](\tau)-E_\eta[u](\tau)}\,d\tau
\,T^{1/{\mathfrak q}} \\
&\leq \norm{v-u}_{\mathfrak p}\norm{E_\eta[v]}_{\mathfrak q} T^{1/{\mathfrak q}}+
\norm{u}_{\mathfrak p}\norm{E_\eta[v]-E_\eta[u]}_{\mathfrak q}\,T^{1/{\mathfrak q}}.
\end{align*}
From the induction hypothesis $E_\eta$ is continuous in the desired sense. Thus, it follows
that for any $\epsilon>0$, there exists a $\delta^\prime_\eta>0$ such that
\begdi
\norm{E_\eta[v]}_{\mathfrak q}\leq \norm{E_\eta[u]}_{\mathfrak q}+1
\text{ and }
\norm{u}_{\mathfrak p}\norm{E_\eta[v]-E_\eta[u]}_{\mathfrak q} T^{1/{\mathfrak q}}<\epsilon/2
\enddi
for all $v$ in a ball centered at $u$ of radius
$\delta^\prime_\eta>0$.\footnote{Of course, $\delta^\prime_\eta$
must be selected so that this ball
is contained inside $\Bpm (R)[0,\,T]$. It is also being tacitly assumed that $u$ is not on the
boundary of $\Bpm (R)[0,\,T]$.
Otherwise, this argument needs a few minor adjustments.}
In which case,
choose
\begdi
\delta_{x_i\eta}=\min\left\{\delta^\prime_\eta,\frac{\epsilon/2}{(\norm{E_\eta[u]}_{\mathfrak q}+1)T^{1/{\mathfrak q}}}\right\}
\enddi
so that if $\norm{u-v}_{\mathfrak p}<\delta_{x_i\eta}$, then
\begdi
\norm{E_{x_i\eta}[v]-E_{x_i\eta}[u]}_{\mathfrak q}<\epsilon.
\enddi
Hence, by induction, $E_\eta$ is continuous
with respect to the $L_{\mathfrak p}$ and
$L_{\mathfrak q}$ norms
for every $\eta\in X^\ast$.

To show that $F_c$ is also continuous in the desired sense,
observe that for any integer $N>0$
\begin{align*}
&\norm{F_c[v]-F_c[u]}_{\mathfrak q}
= \left\Vert\sum_{k=0}^\infty\sum_{\eta\in X^k}
(c,\eta)(E_\eta[v]-E_\eta[u])\right\Vert_{\mathfrak q} \\
\leq& \left\Vert\sum_{k=0}^{N-1}\sum_{\eta\in X^k}
(c,\eta)(E_\eta[v]-E_\eta[u])\right\Vert_{\mathfrak q}+ 
\left\Vert\sum_{k=N}^\infty\sum_{\eta\in X^k}
(c,\eta)(E_\eta[v]-E_\eta[u])\right\Vert_{\mathfrak q} \\
\leq& \left\Vert\sum_{k=0}^{N-1}\sum_{\eta\in X^k}
(c,\eta)(E_\eta[v]-E_\eta[u])\right\Vert_{\mathfrak q}+ 2\sum_{k=N}^\infty K(MR(m+1))^k,
\end{align*}
where $c \in  \ellInfty{M}{X^*}{\R^\ell}$ and $K = \norm{c}_{\ell_{\infty,M}}>0$.
Clearly the second term above can be bounded by $\epsilon/2$ by selecting $N$
to be sufficiently large. Having done this, take $\delta := \min_{|\eta| \leq N} \delta_\eta$, where
the $\delta_\eta$ have been chosen as above to bound the first term by $\epsilon/2$.
This establishes the continuity of the map to $L^\ell_{\mathfrak{q}}[t_0,t_0+T]$. Moreover, in light of \eqref{e-tl2}, it follows that $\norm{F_c[u]}_{\mathfrak{q}} \leq \lVert c \rVert_{\ell_{\infty,M}}\sum_{n\geq 0} (MR(m+1))^n$, where the series is a convergent geometric series. Hence, there exists a constant $S>0$ depending only on $\lVert c\rVert_{\ell_{\infty,M}}$ bounding the $L_\mathfrak{q}$-norm of $F_c[u]$.
\endpr

Now the stronger property of joint continuity is derived using some of concepts developed for the previous two theorems.
First recall that for Banach spaces $V,W$ and $U \subseteq V$ open, the following spaces are Banach spaces:
\begin{itemize}
\item $\text{BC}(U,W)$ the space of bounded continuous functions with the supremum norm $\lVert \cdot \rVert_\infty$.
\item $L(V,W)$ the space of bounded linear functions with the operator norm $\lVert \cdot \rVert_{\text{op}}$.
\end{itemize}

\begth \label{jointBC} (joint continuity)
Let $M \in \re^+$, $\mathfrak p , \mathfrak q$ be conjugate exponents, and $0< T \leq R < 1/(M(m+1))$. The maps
\begin{align*}
 \Phi \colon &\ellInfty{M}{X^*}{\R^\ell} \rightarrow \text{BC}(\Bpm (R)[t_0, \, t_0+T], L_{\mathfrak q}^\ell [t_0,t_0+T]),\quad c \mapsto F_c\\
 \Psi \colon &\Bpm (R)[t_0, \, t_0+T] \rightarrow L(\ellInfty{M}{X^*}{\R^\ell}, L_{\mathfrak q}^\ell [t_0,t_0+T]),\quad  u \mapsto ( c \mapsto F_c[u]      )
\end{align*}
are well defined and continuous. Therefore, the joint map
\begin{align}
\ellInfty{M}{X^*}{\R^\ell}\times \Bpm (R)[t_0, \, t_0+T] &\rightarrow L_{\mathfrak q}^\ell [t_0,t_0+T],\quad (c,u) \mapsto F_c [u] \label{jointmap1}
\end{align}
is also continuous.
\endth

\begpr
First observe from Theorem \ref{th:Fc_operator_continuity} that for every $c \in \ellInfty{M}{X^*}{\R^\ell}$, $F_c$ is bounded and continuous on $\Bpm (R)[t_0, \, t_0+T]$.
Hence, $\Phi$ is well defined and clearly linear.
Furthermore, from \rref{e-tl2},
\begin{align*}
 \lVert \Phi(c)\rVert_\infty = \sup_{u \in \Bpm (R)[t_0, \, t_0+T]} \lVert F_c [u]\rVert_{\mathfrak{q}} \leq \sum_{k=0}^\infty \lVert c\rVert_{\ell_{\infty,M}}(MR(m+1))^k.
\end{align*}
Suitably choosing $R$, the right hand side will be bounded by a finite constant times the factor
$\lVert c \rVert_{\ell_{\infty,M}}$. Hence, $\lVert \Phi \rVert_{\text{op}}< \infty$, and $\Phi$ is continuous.

In light of Theorem~\ref{thm:outputcont}, $\Psi$ is also well defined. To see that $\Psi$ is continuous, let $\epsilon >0$ and fix $u, v \in \Bpm (R)[t_0, \, t_0+T]$.
Then the estimate in the proof of Theorem \ref{th:Fc_operator_continuity} yields
\begin{align*}
 \norm{\Psi (u) - \Psi(v)}_{\text{op}} &= \sup_{\lVert c\rVert_{\ell_{\infty,M}}=1} \norm{F_c [u] - F_c [v]}_{\mathfrak q} \\
 &\leq \sum_{k=0}^{N-1}\sum_{\eta\in X^k} \underbrace{|(c,\eta)|}_{\leq M^{k}k!} \left\Vert E_\eta[v]-E_\eta[u]\right\Vert_{\mathfrak{q}}+  2\sum_{k=N}^\infty (MR(m+1))^k.
\end{align*}
Choosing $N$ sufficiently large, the second term is smaller than $\epsilon/2$. It is known from the proof of Theorem \ref{th:Fc_operator_continuity} that every
$E_\eta$ is continuous, hence one can choose $\delta >0$ such that the first term is less than $\epsilon /2$ if $v$ is in the $\delta$-ball around $u$. Therefore,
$\Psi$ is continuous.

The continuity of the joint map \eqref{jointmap1} follows directly from the continuity of $\Psi$ by \cite[I. \S 3 Proposition 3.10]{MR1666820}. However, it can also be easily derived from the continuity of $\Phi$ as shown next. Select $c_1,c_2 \in \ellInfty{M}{X^*}{\R^\ell}$ and $u,v \in \Bpm (R)[t_0, \, t_0+T]$.
 Applying the triangle inequality, observe that $\norm{F_{c_1}[u]-F_{c_2} [v]}_{\mathfrak q}$ is dominated by
\begin{align}\begin{aligned}
 &\norm{\underbrace{F_{c_1}[u_1]- F_{c_2}[v]}_{= \Phi(c_1-c_2)(u)}}_{\mathfrak q} +\norm{F_{c_1}[u] -F_{c_2}[v]}_{\mathfrak q}  \\
 &\leq \lVert \Phi\rVert_{\text{op}} \lVert c_1-c_2\rVert_{\ell_{\infty,M}}+\norm{F_{c_2}[u] -F_{c_2}[v]}_{\mathfrak q}. \end{aligned}\label{domsum}
\end{align}
Finally, Theorem \ref{th:Fc_operator_continuity} shows that \eqref{domsum} converges to $0$ as $c_1$ tends to $c_2$ and $u$ to $v$, hence \eqref{jointmap1} is continuous.
\endpr

\begex
A nonlinear system identification problem is solved in \cite{Gray-etal_arxiv19,Venkatesh-etal_19} by truncating \rref{eq:Fliess-operator-defined} up to words of length
$J$ and then applying a recursive least-squares algorithm to identify the coefficients of the generating polynomial $p:=\sum_{\eta\in X^{\leq J}}(c,\eta)\eta$.
In this case, Theorem~\ref{jointBC} applies directly to the sequence of estimates $(\hat{p}_j)_j$ with, for example, $M=1$.
\endex

\begex
\cite{Cuchiero-etal_19} describe a model for deep neural networks using \rref{eq:general-MIMO-control-affine-system} with
parameter dependent vector fields $g_i(z,\theta)$. Here $\theta$ is assumed to be the set of fixed parameters of the network, while the inputs $u_i$ correspond to
the \emph{trainable} parameters. In light of \rref{eq:c-equals-Lgh}, if these vector fields are analytic in the state, such networks constitute
a family of Chen-Fliess series $C=\{c_\theta\in\allseriesellLC:\theta\in\Theta\}$. In which case, each generating series would have a $\theta$
dependent growth parameter $M(\theta)$. For a fixed $\theta\in\Theta$, Theorem~\ref{th:Fc_operator_continuity} ensures that the output of the network is a
continuous function of the trainable parameters. In addition, if $\sup_{\theta\in\Theta} M(\theta)$ is finite, then Theorem~\ref{jointBC} guarantees
joint continuity in both the trainable parameters and over the set of design parameters $\Theta$.
\endex

\section{Conclusions}

Sufficient conditions were given under which input-output systems represented by
locally convergent Chen-Fliess series are \emph{jointly} continuous
with respect to their generating series and as operators mapping a ball in an $L_{\mathfrak p}$-space to a ball in an $L_{\mathfrak q}$-space,
where $\mathfrak p$ and $\mathfrak q$ are conjugate exponents.
Continuity with respect to the generating series was characterized using Banach topologies on subsets of $\allseriesell$.
These results were then combined
with elements of proofs for other
forms of continuity appearing
in the literature to produce the desired joint continuity result.

\addcontentsline{toc}{section}{References}
\bibliographystyle{new}
\bibliography{fliess_lit}

\end{document}